\documentclass[conference,letterpaper]{IEEEtran}

\usepackage{cite}
\usepackage{textcomp}

\newenvironment{changemargin}[2]{\begin{list}{}{%
\setlength{\topsep}{0pt}%
\setlength{\leftmargin}{0pt}%
\setlength{\rightmargin}{0pt}%
\setlength{\listparindent}{\parindent}%
\setlength{\itemindent}{\parindent}%
\setlength{\parsep}{0pt plus 1pt}%
\addtolength{\leftmargin}{#1}%
\addtolength{\rightmargin}{#2}%
}\item }{\end{list}}

\makeatletter
\let\NAT@parse\undefined
\makeatother
\usepackage{hyperref}  
\hypersetup{
    colorlinks= true,
    citecolor= black,
    linkcolor= black, 
    pdfborder= {0 0 0},
}

\usepackage{cite}
\usepackage{textcomp}
\usepackage{amsmath, amssymb, amsfonts}
\usepackage{bbm}
\usepackage{accents}
\usepackage{algorithm}
\usepackage{algcompatible}
\usepackage{algpseudocode}
\usepackage{color}
\usepackage{graphicx}
\usepackage{enumerate}
\usepackage{booktabs}
\usepackage{subcaption}
\usepackage{enumerate}
\usepackage{mathrsfs}
\usepackage{mathtools}
\usepackage{dsfont}
\usepackage{relsize}

\makeatletter
\let\NAT@parse\undefined
\makeatother
\usepackage{hyperref}  
\hypersetup{
    colorlinks= true,
    citecolor= black,
    linkcolor= black, 
    pdfborder= {0 0 0},
}

\newtheorem{lemma}{Lemma}

\newtheorem{theorem}{Theorem}
\newtheorem{proposition}{Proposition}

\begin{document}
\title{Semantics of Instability in Networked Control}

\makeatletter
\author{Saad Kriouile$^1$, Mohamad Assaad$^1$, and Touraj Soleymani$^{2}$\\
{\tt\small \{saad.kriouile, mohamad.assaad\}@centralesupelec.fr}, {\tt\small touraj@city.ac.uk}\\[1.5\jot]
$^1$ CentraleSupelec, University of Paris-Saclay, France\\
$^2$ City St George's, University of London, United Kingdom
}

\maketitle

\begin{abstract}
This paper addresses a scheduling problem in the context of a cyber-physical system where a sensor and a controller communicate over an unreliable channel. The sensor observes the state of a source at each time, and according to a scheduling policy determines whether to transmit a compressed sampled state, transmit the uncompressed sampled state, or remain idle. Upon receiving the transmitted information, the controller executes a control action aimed at stabilizing the system, such that the effectiveness of stabilization depends on the quality of the received sensory information. Our primary objective is to derive an optimal scheduling policy that optimizes system performance subject to resource constraints, when the performance is measured by a dual-aspect metric penalizing both the frequency of transitioning to unstable states and the continuous duration of remaining in those states. We formulate this problem as a Markov decision process, and derive an optimal multi-threshold scheduling policy.
\end{abstract}


\section{Introduction}
Advances in communication, computation, and control technologies have driven the development of cyber-physical systems~\cite{kim2012cyber}. These systems often exhibit dynamic and distributed characteristics, necessitating frequent status updates of remote controllers. Such continuous inflow of information enables cyber-physical systems to promptly adjust their behaviors. However, various constraints can impede real-time status updating, creating discrepancies between the system's actual states and its perceived states~\cite{touraj2023tac, touraj2024tit}. This mismatch itself can subsequently place the system in unstable conditions\footnote{The term unstable here can be interpreted also as undesirable}~\cite{tatikonda2004control, ranade2018control}. In this paper, we propose a dual-aspect semantic framework ~\cite{uysal2022semantic} that improves the performance of cyber-physical systems in the above conditions. Our framework, which is appropriate for stabilization tasks over communication channels, penalizes both the frequency of transitioning to unstable states and the continuous duration of remaining in those states. We refer to these aspects shortly as the presence and persistence of instability, i.e.,
\begin{enumerate}
    \item \textbf{Presence of Instability}: This aspect quantifies the instances when the system transitions into an unstable state. Frequent transitions to instability indicate situations where the system tends to revert to undesirable conditions even after corrective mechanisms have been applied.
    \item \textbf{Persistence of Instability}: This aspect accounts for instances where the system remains in an unstable state. Prolonged instability can significantly degrade system performance, underscoring the importance of corrective mechanisms that rapidly resolve anomalies.
\end{enumerate} 

Note that a dual-aspect metric---penalizing both the presence and persistence of instability and their interactions---is highly relevant across various applications. For instance, in autonomous driving systems, minimizing the frequency of transitioning to and the continuous duration of remaining in potentially hazardous conditions is essential for safety and efficiency. Therefore, a metric that takes into account both the presence and persistence aspects can inform better navigational decision-making~\cite{alvaro2016prolonged}, ultimately enhancing overall system performance.

This paper addresses a scheduling problem in the context of a cyber-physical system where a sensor and a controller communicate over an unreliable channel. The sensor observes the state of a source at each time, and according to a scheduling policy determines whether to transmit a compressed sampled state, transmit the uncompressed sampled state, or remain idle. Upon receiving the transmitted information, the controller executes a control action aimed at stabilizing the system, such that the effectiveness of stabilization depends on the quality of the received sensory information. Our primary objective is to derive an optimal scheduling policy that optimizes system performance subject to resource constraints, when the performance is measured by a dual-aspect metric addressing both the presence and persistence of instability. This framework ensures effective stabilization despite communication constraints.

Our setting correlates stabilization effectiveness with data quality, and is motivated by networked systems in which transmitting the entire collected data to the controller is both bandwidth-intensive and energy-consuming. Consider a satellite in the orbit, where onboard sensors collect data about the vehicle and surrounding obstacles. To maintain stable maneuver and avoid collisions, the satellite needs to communicate this information to a ground-based controller that assists in navigation and decision-making. Instead of transmitting raw high-dimensional sensory information such as visual data, the satellite can extract and transmit critical features, such as deviations from a predefined trajectory or proximity to obstacles. Transmitting this compressed low-dimensional sensory information reduces bandwidth and energy usage. However, poor sensory information may limit the controller's ability to perform fine-grained adjustments, potentially affecting the satellite’s performance in complex missions.

\subsection{Related Work}
Several semantic metrics have been developed for networked real-time systems to evaluate system performance in undesirable conditions~\cite{kaul2012age,kriouile2021global,yates2018age,soleymani202211,sun2017age2,sun2017, maatouk2020lexicographic,kriouile2021aoii, kam2020age, yates2021age,maatouk2020age, maatouk2022age, voi,voi2}. The most common state-agnostic metric is the age of information (AoI), introduced in \cite{kaul2012age}, which measures the freshness of information sent from a sensor to a receiver. Specifically, AoI increases whenever the destination does not receive new status updates about the tracked source, irrespective of the source's state. Another state-agnostic metric is the version age of information (VAoI), introduced in \cite{yates2018age}, which quantifies the lag in version updates between the source and the receiver. To address the limitations of AoI, state-aware metrics such as the age of incorrect information (AoII)~\cite{maatouk2020age} have been introduced. AoII considers the relevance of status updates by evaluating the correctness of the information received. Similarly, the version innovative age of information (VIAoI), proposed in \cite{salimnejad2024age}, refines VAoI by measuring the number of outdated versions at the receiver compared to the source when the source is in specific states. In \cite{kriouile2023pull}, an instance of AoII was proposed to account for the distances between different states of a Markov process. However, as noted in \cite{luo2024exploiting}, these metrics are generally content-agnostic, treating all source states equally. This limitation makes them less effective in applications where certain states are more critical than others. To overcome this issue, content-aware metrics have been proposed. For example, \cite{stamatakis2019control} introduces an age of information function that penalizes delays in updating alarm states more heavily. Furthermore, \cite{luo2024exploiting} introduces the age of missed alarm (AoMA) and the age of false alarm (AoFA) metrics, which account for the lasting impacts of missed and false alarm errors, respectively. These works focus primarily on real-time monitoring applications, where the main objective is to track a source rather than control it. Specifically, the actions taken by the controller in these scenarios do not influence the monitored process directly.

In contrast to these works, our paper considers that the central entity not only monitors the source but also has the ability to stabilize it when it is in undesirable condition. To address this, we introduce a novel framework in which the source is sampled by taking into account the presence and persistence of the source in unstable states and the effectiveness of stabilization based on the quality of communicated data. We consider an age-based metric similar to the ones used in \cite{luo2024exploiting,maatouk2022age}), which captures the Age of System Instability (AoSI) in the context of a stabilization task. The main novelty in our framework lies in the fact that the controller can impact the evolution of the source, which makes the underlying analysis more challenging.

\subsection{Paper Outline}
In section \ref{sec:syst_mod}, we describe our system model, then in Section \ref{sec:met_evol}, we introduce the AoSI that captures the presence and persistence of instability, and derive its evolution in a Markov decision process (MDP) framework. In Section \ref{sec:prob_form}, we formulate the problem of interest with the goal to minimize a linear combination of the average AoSI and the average energy consumed by the sensor. In Section \ref{sec:optimal_sol}, we derive the structure of the optimal solution. Finally, we present numerical results in Section \ref{sec:num_reslt}, and conclude the paper in Section \ref{sec:conclusion}.

\begin{figure}[t!]
    \centering
    \includegraphics[width=0.55\linewidth]{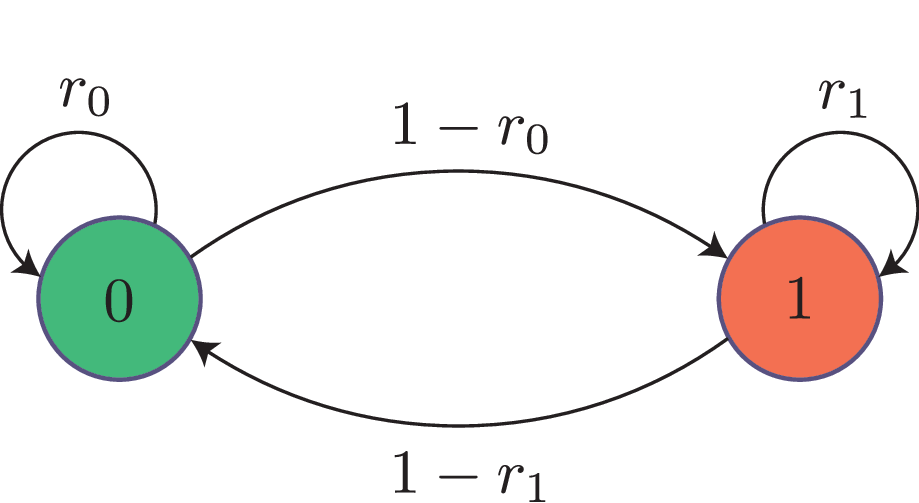}
    \caption{A binary Markov process, where the states represent the macroscopic modes of the source.}
    \label{fig:markov_chain}
\end{figure}

\section{System Model}\label{sec:syst_mod}
\color{black}
\subsection{Networked System Description}
Our networked system consists of a dynamical process (source), a sensor (transmitter), and a controller (receiver). The dynamical process evolves within a high-dimensional state space, which includes both stable (desired) states and unstable (undesired) states. We assume that the probability of transitioning from a stable state to an unstable state is $1-r_0$ where $r_0+r_1 \ge 1$,
and the probability of transitioning from an unstable state to a stable state is $1-r_1$. As such, the dynamical process can be modeled as a binary Markov process, where the states represent the macroscopic modes of the source (see Fig.~\ref{fig:markov_chain}). We assume that, in this Markov process, state $0$ corresponds to the stable mode, while state $1$ corresponds to the unstable mode. The sensor observes the state of the source at each time and decides whether to transmit information about the current state to the controller. The sensor has three actions available at each time: transmitting compressed information, transmitting uncompressed information, or staying idle. We denote by $\delta(t)$ the action taken by the sensor that equals to $1$ if the compressed information about the source state is transmitted, equals to $2$ if the uncompressed information is transmitted, and equals to $0$ if there is no transmission. A transmitted packet is successfully decoded by the controller with a probability $\rho$.

If the source is in the unstable mode, transmitting the uncompressed message is more effective than transmitting the compressed message for stabilizing the source. The uncompressed message provides richer information, enabling the controller to make better decisions to stabilize the source. In our model, we represent this stabilization effectiveness by assigning a higher probability $q$ to the likelihood of the source transitioning to the stable mode when the controller receives an uncompressed message (high-quality information), compared to the probability $p$ when the controller receives a compressed message (low-quality information), i.e., $q > p$. However, this improvement comes at the cost of higher energy consumption. Specifically, transmitting the uncompressed message requires more energy, denoted by $\lambda_2$, compared to the energy required for the compressed message, denoted by $\lambda_1$, i.e., $\lambda_2>\lambda_1$.

The objective in our work is to keep the source in the stable mode as much as possible by penalizing presence and persistence of the unstable mode subject to constraints on communication energy.
To evaluate the performance of our system, we adopt a metric similar to the ones developed in \cite{luo2024exploiting,maatouk2022age}, which we refer here to as the Age of System Instability (AoSI). This metric increases as long as the source is in or remains in an unstable state, and resets to 0 once the system stabilizes. We describe in the following section the detailed evolution of AoSI within our novel  framework where the monitor has control over the source.

\subsection{Semantic Metric Evolution}\label{sec:met_evol}
We describe how AoSI evolves in our system settings. Let $s(t)$ be the value of AoSI at time $t$. We distinguish between two cases of $s(t)$, $s(t)=0$ and $s(t)>0$. When $s(t)=0$, the source is in the stable mode. Therefore, if the sensor decides to stay idle, AoSI will be incremented by one if the source moves to the other mode. If the sensor decides to transmit a status update, AoSI will be incremented by one in the two following scenarios: if the packet is unsuccessfully decoded by the controller and the source moves to the unstable mode; or if the packet is successfully decoded, the source moves to the other mode but the controller does not succeed to stabilize the source to the stable mode.

When $s(t)\neq 0$, the source is in the unstable mode. Hence, if the sensor decides to stay idle, AoSI will be incremented by one if the source remains at the unstable mode. If the sensor decides to transmit a status update, AoSI will be incremented by one in the two following scenarios: if the packet is unsuccessfully decoded by the controller and the source remains at the same mode; or if the packet is successfully decoded, the source remains at the same mode but the controller does not succeed to stabilize the source to the stable mode.

Given that AoSI at time $t$ can change only to two different values, $s(t)+1$ and $0$, the sum of the two transition probabilities; i.e.; the transition probabilities to $s(t)+1$ and to $0$, is equal to $1$. Thus, for sake of space, we present in the following only the transition probabilities to the state $s(t)+1$:
\begin{itemize}
    \item $\Pr\big(s(t+1)=s(t)+1|s(t)=0,\delta(t)=0\big)=1-r_0$.
    \item $\Pr\big(s(t+1)=s(t)+1|s(t)=0,\delta(t)=1\big)=(1-\rho)(1-r_0)+\rho (1-p) (1-r_0)$.
    \item $\Pr\big(s(t+1)=s(t)+1|s(t)=0,\delta(t)=2\big)=(1-\rho)(1-r_0)+\rho (1-q) (1-r_0)$.
    \item $\Pr\big(s(t+1)=s(t)+1|s(t)>0,\delta(t)=0\big)=r_1$.
    \item $\Pr\big(s(t+1)=s(t)+1|s(t)>0,\delta(t)=1\big)=(1-\rho)r_1+\rho (1-p) r_1$.
    \item $\Pr\big(s(t+1)=s(t)+1|s(t)>0,\delta(t)=2\big)=(1-\rho)r_1+\rho (1-q) r_1$.
\end{itemize}
Other probabilities can be deduced given the fact that $\Pr\big(s(t+1)=0|s(t),\delta\big)=1-\Pr\big(s(t+1)=s(t)+1|s(t),\delta\big)$.

\subsection{Problem Statement}\label{sec:prob_form}
Our goal is to determine the optimal scheduling policy, that minimizes a linear combination of the average AoSI and the average energy usage. A scheduling policy is defined as $\pi = (\delta^{\pi}(0), \delta^{\pi}(1), \dots)$, where $\delta^{\pi}(t)$ refers to the action taken by the sensor at time $t$. Taking into account the energy cost incurred when sending a packet, the stage cost at time $t$ is $C(t) = s(t) + \lambda_1 \mathbbm{1}_{\delta(t)=1}+ \lambda_2 \mathbbm{1}_{\delta(t)=2}$. Mathematically, our objective is to solve the following optimization problem:
\begin{align}\label{eq:problem_formulation}
& \underset{\pi\in \Pi}{\text{minimize}} \ \lim_{T\to\infty} \text{inf}\:\frac{1}{T}\mathbb{E}^{\pi\in \Pi}\Big[ \sum_{t=0}^{T-1}C^{\pi}(t) \big|s(0)\Big]
\end{align}
where $\Pi$ denotes the set of all causal scheduling policies and $C^{\pi}(t)$ denotes the stage cost function at time $t$ under the scheduling policy~$\pi$.

\section{Main Result}\label{sec:optimal_sol}
In this section, we derive the structure of the optimal solution. Building on the procedure used in \cite{kriouile2021aoii,maatouk2020age,kriouile2023pull,kriouile2021global,larranaga2017asymptotically}, we formulate our problem as an MDP. Then, leveraging a backward induction, we derive the optimal solution.

\subsection{Reformulation to an MDP and Optimal Solution Structure}
The optimization problem in (\ref{eq:problem_formulation}) can be viewed as an infinite-horizon average-cost MDP problem with the following characteristics:
\begin{itemize}
\item \emph{State}: The state at time $t$, denoted by $s(t)$, represents AoSI, depending on the context. 
\item \emph{Actions}: The actions at time $t$, denoted by $\delta(t)$ representing the transmission decision.
\item \emph{Transition}: The transition probabilities specify the probabilities associated with state changes given actions.
\item \emph{Cost}: The instantaneous cost, denoted by $C(s(t),\delta(t))$, is equal to $s(t)+\lambda_1 \mathbbm{1}_{\delta(t)=1}+\lambda_2 \mathbbm{1}_{\delta(t)=2}$.
\end{itemize}
The optimal scheduling policy $\pi^*$ of the problem in \eqref{eq:problem_formulation} can then be obtained by solving the following Bellman equation for each state $s$:
\begin{align}
\theta& + V(s) \nonumber = \min_{\delta \in\{0,1,2\}} \Big\{s+\lambda_1 \mathbbm{1}_{\delta=1}+\lambda_2 \mathbbm{1}_{\delta=2}\\
&\qquad \qquad \qquad \qquad \qquad +\sum_{s'\in \mathbb{N} }\Pr(s\rightarrow s'|\delta)V(s') \Big\}
\label{eq:bellman_general}
\end{align}
where $\Pr(s\rightarrow s'|\delta)$ is the transition probability from $s$ to $s'$ under action $a$, $\theta$ is the optimal value, $V(s)$ is the differential cost-to-go function when the state is $s$, and $\mathbb{N}$ is the set of positive integers. The next Lemma and Theorem, which are the main results of this paper, specify the structure of the optimal scheduling~policy.
\begin{lemma}\label{lem:V_increasing}
$V(s)$ is an increasing function with respect to $s$.
\end{lemma}
\begin{IEEEproof}
See Appendix \ref{app:lem:V_increasing}.
\end{IEEEproof}
\begin{theorem}\label{theo:threshold_policy}
The optimal solution to the problem in (\ref{eq:bellman_general}) is an increasing multi-threshold policy. There exists $n_1$ and $n_2$ with $n_1 \le n_2$ such that the optimal action is $\delta=0$ if $s< n_1$, the optimal action is $\delta=1$ if $n_1 \le s <n_2$, and the optimal action is $\delta=2$ for if $s \ge n_2$.
\end{theorem}
\begin{IEEEproof}
See Appendix \ref{app:theo:threshold_policy}
\end{IEEEproof}

\subsection{Computation of Optimal Multi-Threshold Policy}
\subsubsection{Closed form expression of the objective function}
After establishing that the structure of the optimal solution is multi-threshold policy, we study now the steady-state form of the optimization problem \eqref{eq:problem_formulation} under a given multi-threshold policy. For that, we first express the steady-state form of the problem in \eqref{eq:problem_formulation} as follows:

\begin{equation}\label{eq:steady_state_form}
\begin{aligned}
& \underset{n\in \mathbb{N}}{\text{minimize}} 
& & \overline{s^{n}}+\lambda_1\overline{\delta_1^n}+ \lambda_2 \overline{\delta^n_2}
\end{aligned}
\end{equation}
where $\overline{s^{n}}$ is the average value of AoSI, $\overline{\delta_1^n}$ and $\overline{\delta^n_2}$ are the average time during which the sensor sends a compressed message and an uncompressed message, respectively, under a multi-threshold policy with parameter~$n=(n_1,n_2)$. More specifically, we have
\begin{align}
\overline{s^{n}}&=\lim_{T\to+\infty} \text{sup}\:\frac{1}{T}\mathbb{E}^{n}\Big(\sum_{t=0}^{T-1}s(t)|s(0),\pi(n)\Big)\label{eq:average_age}\\
\overline{\delta_1^n}&=\lim_{T\to+\infty} \text{sup}\:\frac{1}{T}\mathbb{E}^{n}\Big(\sum_{t=1}^{T}\mathbbm{1}_{\delta(t)=1}|s(0),\pi(n)\Big)\label{eq:average_active_time}\\
\overline{\delta_2^n}&=\lim_{T\to+\infty} \text{sup}\:\frac{1}{T}\mathbb{E}^{n}\Big(\sum_{t=1}^{T}\mathbbm{1}_{\delta(t)=2}|s(0),\pi(n)\Big)
\end{align}
where $\pi(n)$ denotes the multi-threshold policy with parameter $n=(n_1,n_2)$.

\begin{figure}[t!]
\centering
\includegraphics[trim={0.4cm -0.5cm 0 0.25cm}, clip, width=1\columnwidth]{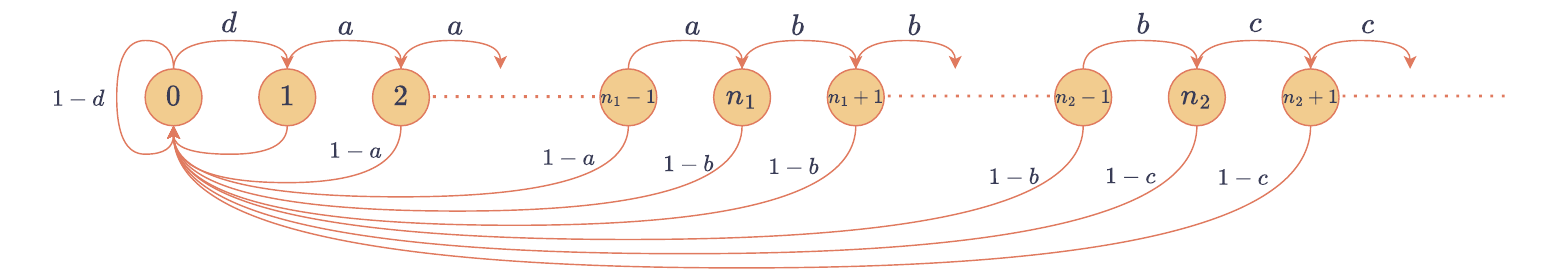}
\caption{The AoSI state transition under a multi-threshold policy with parameter $n=(n_1,n_2)$, when $n_1, n_2>0 $.}
\label{fig:dtmc_n1_n2_diff_0}
\end{figure}

\begin{figure}[t!]
\centering
\includegraphics[trim={0 0.5cm 0 0.25cm}, clip, width=0.9\columnwidth]{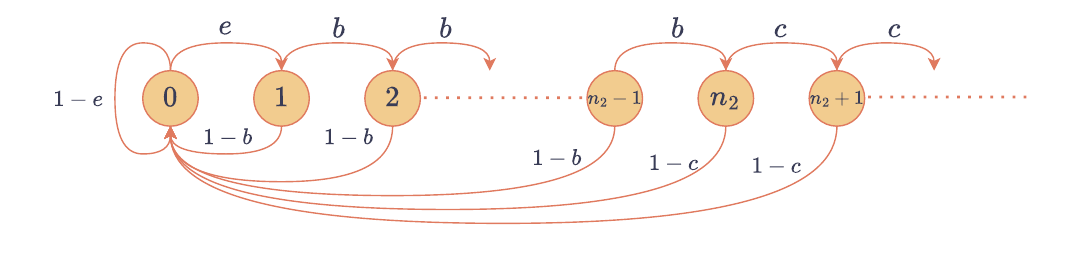}
\caption{The AoSI state transition under a multi-threshold policy with parameter $n=(n_1,n_2)$, when $n_1=0$ and $n_2>0$.}
\label{fig:dtmc_n1_eq_0_n2_diff_0}
\end{figure}

\begin{figure}[t!]
\centering
\includegraphics[trim={0 0 0 0.25cm}, clip, width=0.65\columnwidth]{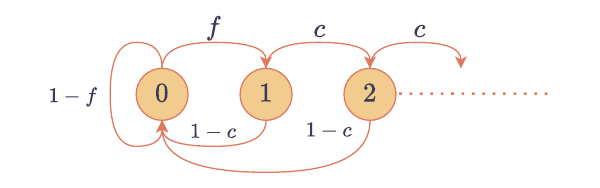}
\caption{The AoSI state transition under a multi-threshold policy with parameter $n=(n_1,n_2)$, when $n_1=0$ and $n_2=0$.}
\label{fig:dtmc_n1_n2_eq_0}
\end{figure}

To compute $\overline{s^{n}}$, $\overline{\delta_1^n}$ and $\overline{\delta_2^n}$, we will first establish the stationary distribution of the discrete-time Markov chain (DTMC) that illustrates the evolution of AoSI under a multi-threshold policy with parameter $n=(n_1,n_2)$ (see Fig.~\ref{fig:dtmc_n1_n2_diff_0}, \ref{fig:dtmc_n1_eq_0_n2_diff_0} and \ref{fig:dtmc_n1_n2_eq_0}, where $a=r_1$, $b=r_1(1-\rho p)$, $c=r_1(1-\rho q)$, $d=1-r_0$, $e=(1-r_0)(1-\rho p)$ and $f=(1-r_0)(1-\rho)$). In he following proposition, we derive the stationary distribution of the our DTMC of interest.
\begin{proposition}\label{prop:stationary_distribution}
For a given threshold pair $n=(n_1,n_2)$, the DTMC admits $u_n(i)$ as its stationary distribution distinguishing between three cases:
\begin{changemargin}{-5 pt}{0 pt}
\begin{itemize}
\item If $n_1,n_2 >0$
\begin{equation}
 u_n(i)=\left\{
    \begin{array}{ll}
        u_n(0) & \text{if} \ i=0  \\
        da^{i-1}u_n(0)& \text{if} \ 1 \leq i \leq n_1  \\
        da^{n_1-1}b^{i-n_1}u_n(0)& \text{if} \ n_1 < i \le n_2\\
        da^{n_1-1}b^{n_2-n_1}c^{i-n_2}u_n(0)& \text{if} \ i > n_2\\
    \end{array}
\right.
\end{equation}
\label{stationarydistribution}
where
$u_n(0)\\
=\frac{(1-a)(1-b)(1-c)}{(1-a+d)(1-c)(1-b)+(1-c)(b-a)d a^{n_1-1}+d(1-a)(c-b)b^{n_2-n_1}a^{n_1-1}}$
\item If $n_1=0$ and $n_1>0$:
\begin{equation}
 u_n(i)=\left\{
    \begin{array}{ll}
        u_n(0) & \text{if} \ i=0  \\
        eb^{i-1}u_n(0)& \text{if} \ 1 \leq i \leq n_2  \\
        eb^{n_2-1}c^{i-n_2}u_n(0)& \text{if} \ n_2 < i\\
    \end{array}
\right.
\end{equation}
\label{stationarydistribution}
where
$u_n(0)=\frac{(1-b)(1-c)}{(1-b+e)(1-c)+(c-b)e b^{n_2-1}}$.
\item If $n_0=n_1=0$:
\begin{equation}
 u_n(i)=\left\{
    \begin{array}{ll}
        u_n(0) & \text{if} \ i=0  \\
        fc^{i-1}u_n(0)& \text{if} \ 1 \leq i \\
    \end{array}
\right.
\end{equation}
\label{stationarydistribution}
where
$u_n(0)=\frac{1-c}{1-c+f}$.

\end{itemize}
\end{changemargin}
\end{proposition}
\begin{IEEEproof}
The proof is based mainly on resolving the full balance equation: $u_n(i)=\sum_{j=0}^{+\infty} \Pr(j \rightarrow i) u_n(j)$. After resolving this equation and given that $\sum_{j=0}^{+\infty} u_n(j)=1$, we get our desired result.
\end{IEEEproof}
Leveraging these results, we provide the closed form expression of $\overline{s^n}+\lambda_1 \overline{\delta_1^n} +\lambda_2 \overline{\delta_2^n}$. We can write
\begin{align}\label{eq:explicit_exp_st_st}
    \overline{s^n}+\lambda_1 \overline{\delta_1^n} +\lambda_2 \overline{\delta_2^n}
    =& \sum_{i=0}^{+\infty} i u_n(i)+ \lambda_1 \sum_{i=n_1}^{n_2-1} u_n(i)
    +\lambda_2\sum_{i=n_2}^{+\infty} u_n(i)
\end{align}
Now, we compute the first term of the previous equation, $\sum_{i=0}^{+\infty} i u_n(i)$. For that, we distinguish between three cases: $n_1, n_2 > 0$; $n_1=0$ and $n_2>0$;  $n_1=n_2=0$. 
\begin{itemize}
\item $n_1, n_2 > 0$:
\begin{align}\label{eq:average_aosi_1}
&\overline{s^{n}}=u_n(0)\frac{d[a^{n_1+1}n_1-n_1 a^{n_1}-a^{n_1}+1]}{(1-a)^2} \nonumber \\
&+u_n(0)\frac{d b a^{n_1-1}[b^{n_2-n_1}(bn_2-n_2-1)-b n_1+n_1+1]}{(1-b)^2} \nonumber \\
&+u_n(0)\frac{dc b^{n_2-n_1} a^{n_1-1}(-cn_2+n_2+1)}{(1-c)^2}
\end{align}
\item $n_1=0$ and $n_2>0$:
\begin{align}\label{eq:average_aosi_2}
\overline{s^{n}}=&u_n(0)\frac{e[b^{n_2+1}n_2-n_2 b^{n_2}-b^{n_2}+1]}{(1-b)^2} \nonumber \\
&+u_n(0)\frac{ecb^{n_2-1}[-cn_2+n_2+1]}{(1-c)^2}
\end{align}
\item $n_1=n_2=0$:
\begin{align}\label{eq:average_aosi_3}
    \overline{s^{n}}=u_n(0)\frac{f}{(1-c)^2}
\end{align}
\end{itemize}
To compute the remaining terms of \eqref{eq:explicit_exp_st_st}, we again distinguish between three cases:
\begin{itemize}
    \item $n_1,n_2>0$:
    \begin{align}\label{eq:average_trans_time_1}
        &\sum_{i=n_1}^{n_2-1}u_n(i)= u_n(0) d a^{n_1-1} \frac{1-b^{n_2-n_1}}{1-b} \nonumber \\
        &\sum_{i=n_2}^{+\infty}u_n(i)= u_n(0) d a^{n_1-1} b^{n_2-n_1} \frac{1}{1-c}
    \end{align}
    \item $n_1=0,n_2>0$:
    \begin{align}\label{eq:average_trans_time_2}
        &\sum_{i=n_1}^{n_2-1}u_n(i)= u_n(0) (1+ e\frac{1-b^{n_2-1}}{1-b}) \nonumber \\
        &\sum_{i=n_2}^{+\infty}u_n(i)= u_n(0) e b^{n_2-1} \frac{1}{1-c}
    \end{align}
    \item $n_1=n_2=0$
    \begin{align}\label{eq:average_trans_time_3}
        \sum_{i=n_1}^{n_2-1}u_n(i)= 0, \hspace{1cm}
        \sum_{i=n_2}^{+\infty}u_n(i)= 1
    \end{align}
\end{itemize}
where the expression of $u_n(0)$ is given in Proposition \ref{prop:stationary_distribution}. Therefore, combining \eqref{eq:average_aosi_1}, \eqref{eq:average_aosi_2}, \eqref{eq:average_aosi_3}, \eqref{eq:average_trans_time_1}, \eqref{eq:average_trans_time_2} and \eqref{eq:average_trans_time_3}, we obtain a closed form expression of the objective function in problem \eqref{eq:problem_formulation}, denoted by $F(n_1,n_2)$ given in \eqref{eq:exp_obj_funct}. Hence, our problem is equivalent to:
\begin{equation}\label{eq:compact_prob_formul}
    \underset{n_1,n_2}{\min}F(n_1,n_2)
\end{equation}
To find $n_1^*$ and $n_2^*$ that minimize the function given in \eqref{eq:compact_prob_formul}, since we have only two dimensional thresholds space, we can adopt an algorithm based on gradient descent.

\begin{table*}
\caption{Expression of $F(n_1,n_2)$}
\centering
\begin{minipage}{0.75\textwidth}
\begin{itemize}
\item $n_1,n_2>0$
\begin{align}
F(n_1,n_2)=&\frac{1}{(1-a+d)(1-c)(1-b)+(1-c)(b-a)d a^{n_1-1}+d(1-a)(c-b)b^{n_2-n_1}a^{n_1-1}}  \bigg [ \frac{(1-b)(1-c)d[a^{n_1+1}n_1-n_1 a^{n_1}-a^{n_1}+1]}{(1-a)} \nonumber \\
&+\frac{d b (1-a)(1-c) a^{n_1-1}[b^{n_2-n_1}(bn_2-n_2-1)-b n_1+n_1+1]}{1-b} \nonumber \\
&+\frac{dc(1-a)(1-b) b^{n_2-n_1} a^{n_1-1}(-cn_2+n_2+1)}{(1-c)} \nonumber \\
&+ \lambda_1 d (1-a)(1-c) a^{n_1-1} (1-b^{n_2-n_1}) + \lambda_2 d a^{n_1-1}(1-a)(1-b) b^{n_2-n_1} \bigg] \nonumber
\end{align}
\item $n_1=0$, $n_2>0$:
\begin{align}
F(n_1,n_2)=&\frac{1}{(1-b+e)(1-c)+(c-b)e b^{n_2-1}} \bigg [ \frac{(1-c)e[b^{n_2+1}n_2-n_2 b^{n_2}-b^{n_2}+1]}{1-b}+\frac{(1-b)ecb^{n_2-1}[-cn_2+n_2+1]}{1-c}\nonumber\\
&+\lambda_1 \big((1-b)(1-c)+ e(1-c)(1-b^{n_2-1})\big) +\lambda_2 (1-b)e b^{n_2-1}\bigg] \nonumber
\end{align}
\item $n_1=n_2=0$
\begin{align}\label{eq:exp_obj_funct}
    F(n_1,n_2)=\frac{f}{(1-c+f)(1-c)}+\lambda_2
\end{align}
\end{itemize}
\medskip
\hrule
\end{minipage}
\end{table*}

\section{Numerical Results}\label{sec:num_reslt}
In this section, we analyze the optimal cost in \eqref{eq:compact_prob_formul} associated with the optimal multi-threshold policy as a function of $\lambda_1$ and $\lambda_2$. Recall that $\lambda_1$ and $\lambda_2$  denote the energy costs incurred by the transmitter when sending low-quality information and high-quality information, respectively.
We suppose that $\lambda_1$ and $\lambda_2$ vary from $0$ to $9$ with step size $1$. We consider the following parameters $r_1=0.9$,
$r_0=0.1$, $\rho=0.1$, $p=1/2$ and $q=0.9$. 

\begin{figure}[h!]
    \centering
    \begin{subfigure}[b]{0.49\linewidth}
        \centering
        \includegraphics[trim={0 0 0 1cm}, clip,width=\linewidth]{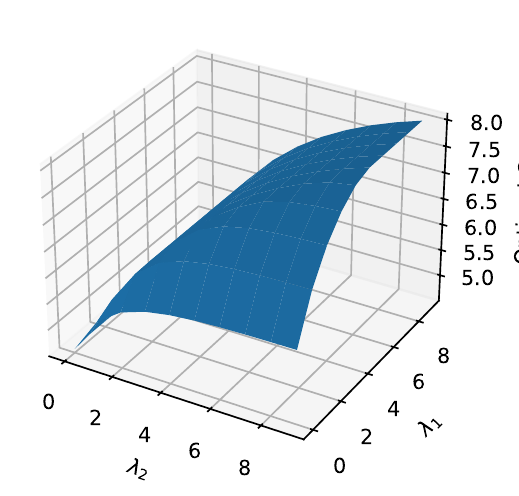}
        \caption{}
        \label{fig:sub1}
    \end{subfigure}
    \hfill
    \begin{subfigure}[b]{0.49\linewidth}
        \centering
        \includegraphics[trim={0 0 0 1cm}, clip,width=\linewidth]{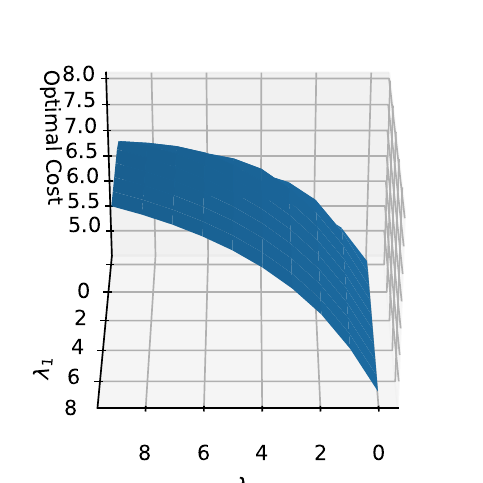} 
        \caption{}
        \label{fig:sub2}    
    \end{subfigure}
    \caption{Optimal Cost in function of $\lambda_1$ and $\lambda_2$}
    \label{fig:evol_opt_cost}
\end{figure}

In Fig.~\ref{fig:evol_opt_cost}, we observe that the average optimal cost increases as $\lambda_1$ and $\lambda_2$ grow. This is because the performance of the optimal solution is reducing as the energy cost of transmitting low-quality or high-quality information increases. Additionally, Fig.~\ref{fig:evol_opt_cost}(b) shows that for $\lambda_2 = 0$, the average cost remains the same regardless of the value of $\lambda_1$. Indeed, when $\lambda_2$ is close to $0$, the best strategy is to always send high-quality information, as the transmission energy cost is negligible, and the superior quality of the information received by the controller enables better decision-making for stabilizing the source. Conversely, when $\lambda_1 = 0$, the optimal solution is not necessarily to always send low-quality information. While the transmission energy cost is negligible, low-quality information can impair the controller's ability to make effective stabilization decisions, leading to suboptimal performance. This explains why, for $\lambda_2 = 0$, the optimal cost remains unchanged (indicating the same optimal policy). However, for $\lambda_1 = 0$, the optimal cost increases with $\lambda_2$, as the reliance on high-quality information comes at a greater energy cost.

\begin{figure}[h!]
\centering
\includegraphics[trim={0 0.25cm 0 0.25cm}, clip, width=0.8\columnwidth]{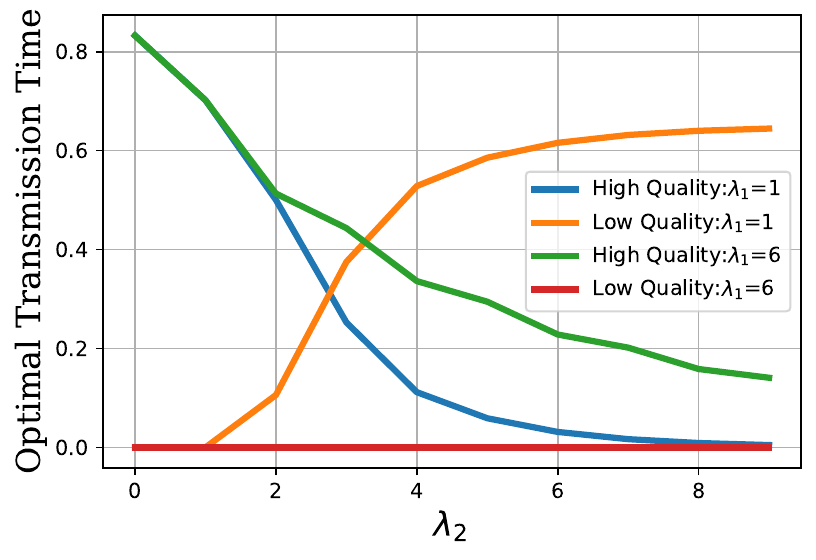}
\caption{The evolution of the average transmission time in function of $\lambda_1$ and $\lambda_2$.}
\label{fig:trans_time_funct_lams}
\end{figure}

In Fig.~\ref{fig:trans_time_funct_lams}, we compare the optimal average transmission time of low-quality information with that of the high-quality information for different values of $\lambda_1$ and $\lambda_2$. Specifically, we plot the curves representing the functions $\sum_{i=n_1^*}^{n_2^*-1} u_{n^*}(i)$ and $\sum_{i=n_2^*}^{+\infty} u_{n^*}(i)$ under the optimal solution $n^*$ in function of $\lambda_2$ for two different values of $\lambda_1$; $1$ and $6$. We deduce that when the energy cost of low-quality information transmission is low ($\lambda_1 = 1$), sending the compressed message becomes more advantageous than sending the uncompressed message as $\lambda_2$ increases, in terms of optimizing the average cost. In contrast, when the energy cost of low-quality information transmission is high ($\lambda_1 = 6$), sending the compressed message is never a favorable strategy. Consequently, the optimal solution in this case involves only two actions: staying idle or transmitting the high-quality information.

\section{Conclusion}\label{sec:conclusion}
In this work, we have investigated a scheduling problem in the context of a cyber-physical system in which sensory information of a source needs to be transmitted from a sensor to a remote controller that is responsible for stabilizing the source. We derived the optimal scheduling policy, striking a balance between performance and constraints. Through rigorous mathematical analysis, we demonstrated that the optimal solution follows a multi-threshold increasing policy. The numerical results further validate our theoretical findings, offering valuable insights into the relationship between the optimal average transmission time for low-quality and high-quality information and the associated transmission energy costs. These results highlight the trade-offs inherent in designing efficient and effective scheduling~policies.

\bibliography{AoS.bib}
\bibliographystyle{ieeetr}

\begin{appendices}
\section{Proof of Lemma \ref{lem:V_increasing}}\label{app:lem:V_increasing}
The relative value iteration equation consists of updating the value function $V_t(.)$ as follows: 
\begin{align}
V_{t+1}(s)=\min\big\{V_t^{0}(s),V_t^{1}(s),V_t^{2}(s)\big\} 
\end{align}
where
\begin{align}
V_t^{0}(s)& = s+r V_t(s+1)+(1-r) V_t(0), \\[1.5\jot]
V_t^{1}(s)& = s+\lambda_1+(1-\rho p)r V_t(s+1)\nonumber \\
& \qquad +((1-\rho p)(1-r)+\rho p) V_t(0), \nonumber\\[1.5\jot]
V_t^{2}(s)& = s+\lambda_2+(1-\rho q)r V_t(s+1)\nonumber \\
& \qquad +((1-\rho q)(1-r)+\rho q) V_t(0), \nonumber\\[1.5\jot]    
\end{align}
where $r=(1-r_0) \mathbf{1}_{s=0}+r_1 \mathbf{1}_{s>0}$.\\
We prove this lemma by induction. We show that $V_t(\cdot)$ is increasing for all $t$ and we conclude that for $V(\cdot)$. As $V_0(.)=0$, then the property holds for $t=0$. If $V_t(.)$ is increasing with respect to $s$, we show that, for $s$, we have $V_{t+1}^{\delta}(s) \leq V_{t+1}^{\delta}(s+1)$, for any $\delta \in \{0,1,2\}$.
We have for $s>0$
\begin{align}
&V_{t+1}^{\delta}(s+1)-V_{t+1}^{\delta}(s) \nonumber\\[1.5\jot]
&\qquad \quad = 1+r'(V_{t}(s+1) - V_{t}(s)).
\end{align}
where $r'=r_1 \mathbf{1}_{\delta=0}+(1-\rho p)r_1 \mathbf{1}_{\delta=1}+(1-\rho q)r_1 \mathbf{1}_{\delta=2}$
Given that $V_t(.)$ is increasing with respect to $s$ and $r'\ge 0$, therefore $V_{t+1}^{\delta}(s+1) - V_{t+1}^{\delta}(s) \geq 0$. As consequence, $V_{t+1}^{\delta}(\cdot)$ is increasing with respect to $s$ when $s>0$. 
Now, we compare $V_{t+1}^{\delta}(1)$ and $V_{t+1}^{\delta}(0)$.
We have:
\begin{align}
V_{t+1}^{0}(1)-V_{t+1}^{0}(0)
=& 1+r_1 V_{t}(2) - (1-r_0)V_{t}(1)\nonumber \\[1.5\jot]
&+(1-r_1-r_0)V_t(0). \nonumber \\[1.5\jot]
\overset{(a)}\ge&  1+(r_1+r_0-1)V_t(1)\nonumber\\
&-(r_0+r_1-1)V_t(0) \nonumber \\
= & 1+(r_0+r_1-1)(V_t(1)-V_t(0)) \nonumber \\
\overset{(b)}{\ge} &  0
\end{align}

where $(a)$ comes from the fact that $V_t(1)\le V_t(2)$ and $(b)$ comes from the fact that $r_0+r_1\ge 1$.
Following the same steps, we show that $V_{t+1}^{\delta}(1)-V_{t+1}^{\delta}(0) \ge 0$ for $\delta=1,2$.

Now, since $V_{t+1}(.)=\min\{V^{0}_{t+1}(\cdot),V^{1}_{t+1}(\cdot),V^{2}_{t+1}(\cdot)\}$, then $V_{t+1}(.)$ is increasing with respect to $s$. We demonstrated by induction that $V_t(.)$ is increasing for all $t$. Knowing that $\lim_{t \to \infty} V_t(s)=V(s)$, we conclude that $V(.)$ is also increasing with respect to $s$.

\section{Proof of Theorem \ref{theo:threshold_policy}}\label{app:theo:threshold_policy}
We compare $V^{0}(s)$ with $V^{1}(s)$ and $V^{0}(s)$ with $V^{2}(s)$.
We have: 
\begin{align}
    V^{1}(s)-V^{0}(s)=\lambda_1+\rho p r(V(0)-V(s+1))
\end{align}
According to Lemma \ref{lem:V_increasing}, this difference is decreasing with $s$ since $V(\cdot)$ is increasing with $s$. That means that there exists $s_1$ such that for all $s < s_1$, $V^{0}(s) < V^{1}(s)$ and for all $s \ge s_1$, $V^{0}(s) \ge  V^{1}(s)$.

Similarly, we have: 
\begin{align}
    V^{2}(s)-V^{0}(s)=\lambda_2+\rho q r(V(0)-V(s+1))
\end{align}
According to lemma \ref{lem:V_increasing}, this difference is as well decreasing with $s$. That means that there exists $s_2$ such that for all $s < s_2$, $V^{0}(s) < V^{2}(s)$ and for all $s \ge s_2$, $V^{0}(s) \ge  V^{2}(s)$.

Therefore, for all $s < \min\{s_1,s_2\}$, the optimal action is $0$ as $\min\{V^{0}(s),V^{1}(s),V^{2}(s)\}=V^{0}(s)$ and for all $s \ge \min\{s_1,s_2\}$, the optimal solution is either $1$ or $2$. Now we prove that for $s\ge \min\{s_1,s_2\}$, the optimal solution is threshold based policy. For that, we compare $V^{1}(s)$ with $V^{2}(s)$:
\begin{align}
    V^{2}(s)-V^{1}(s)=\lambda_2-\lambda_1+\rho (q-p) r(V(0)-V(s+1))
\end{align}
According to lemma \ref{lem:V_increasing}, this difference is as well decreasing with $s$ since $q\ge p$. That means that there exists $s_3$ such that for all $s < s_3$, $V^{1}(s) < V^{2}(s)$ and for all $s \ge s_3$, $V^{1}(s) \ge  V^{2}(s)$. Therefore, for $s \ge \max\{s_3,\min\{s_1,s_2\}\}$, the optimal solution is $2$ and for $\min\{s_1,s_2\} \le s < \max\{s_3,\min\{s_1,s_2\}\}$, the optimal solution is $1$. As consequence, for $n_1=\min\{s_1,s_2\}$ and for $n_2=\max\{s_3,\min\{s_1,s_2\}\}$, we have for all $s< n_1$, the optimal action is $\delta=0$; for all state $n_1 \le s <n_2$, the optimal action is $\delta=1$ and for all $s \ge n_2$, the optimal action is $\delta=2$. That concludes the proof.
\end{appendices}
\end{document}